\documentclass[pdftex,fleqn,10pt]{wlscirep}
\pdfoutput=1
\usepackage{cleveref}
\usepackage[percent]{overpic}
\usepackage[utf8]{inputenc}
\usepackage[T1]{fontenc}
\usepackage{longtable}
\usepackage{comment}
\usepackage[colorinlistoftodos]{todonotes}[1]

\title{Can diffusion alone explain brain-wide distribution of a CSF tracer within 24 hours?}

\author[1]{Lars Magnus Valnes}
\author[2]{Sebastian K. Mitusch}
\author[3]{Geir Ringstad}
\author[4,5]{Per Kristian Eide}
\author[2]{Simon W. Funke}
\author[1,2,*]{Kent-Andre Mardal}

\affil[1]{Department of Mathematics, University of Oslo, Norway }
\affil[2]{Center for Biomedical Computing, Simula Research Laboratory, Lysaker, Norway}
\affil[3]{Division of Radiology and Nuclear Medicine, Department of Radiology - Rikshospitalet, Oslo, Norway}
\affil[4]{Institute of Clinical Medicine, Faculty of Medicine, University of Oslo, Oslo, Norway}
\affil[5]{Department of Neurosurgery, Oslo University Hospital – Rikshospitalet, Oslo, Norway}
\affil[*]{kent-and@math.uio.no}

\begin{abstract}
The recently proposed glymphatic system suggests that bulk flow is important for clearing waste from the brain, and as such
may underlie the development of e.g. Alzheimer's disease. The glymphatic hypothesis is still controversial and several biomechanical modeling studies at the micro-level have at least partially dismissed the system and its assumptions. In contrast, at the macro-level, there are many 
experimental findings in support of bulk flow. Here, we will investigate to what extent the CSF tracer distributions seen in novel magnetic resonance imaging (MRI) investigations over hours and days
are suggestive of bulk flow or diffusion. 
In order to include the complex geometry of the brain, the heterogeneous CSF flow around the brain, and the transport over the 
time-scale of days, we employed the methods of partial differential constrained optimization to identify the apparent diffusion coefficient (ADC) that would correspond best to the MRI findings.  We found that the computed ADC in grey and white matter was respectively 23\% and 82\% larger than the ADC estimated with DTI, which suggests that diffusion may not be the only mechanism governing transport. 
\end{abstract}
\begin{document}

\flushbottom
\maketitle
%
%
\thispagestyle{empty}

\section*{Introduction}


Most types of dementia are associated with accumulation of metabolic by-products within the brain. 
In contrast to the rest of the body, the brain lacks a lymphatic system to clear  these by-products. 
In 2012, a new pathway, called the paravascular pathway, was proposed~\cite{iliff2012paravascular}, which enables efficient brain-wide circulation and clearance. The network of paravascular pathways in the brain was named the glymphatic system as it resembles the lymphatic system in the rest of the body, while the 'g' in glymphatic highlights the importance of the supportive glia cells in the brain.
The paravascular pathways consists of cerebrospinal fluid flowing in parallel with the vasculature in paravascular spaces.
These pathways have the potential to facilitate  exchange between the cerebrospinal fluid (CSF) and the extracellular fluid deep within the brain.  

To what extent and at what scale the glymphatic system accelerates transport compared to extracellular diffusion is still controversial, and several computational modeling studies have dismissed parts of the system at  micro-scale. For example, the previous studies~\cite{holter2017interstitial, smith2017glymphatic} suggest that diffusion dominates in the interstitium. Furthermore, ~\cite{asgari2016glymphatic, brynjfm, Diem} have found that dispersion in the paravascular spaces adds less than a factor two to diffusion for solute transportation. 
However, multiple experimental and imaging findings at the micro-level point towards transport being different and faster than diffusion~\cite{iliff2012paravascular, mestre2018flow, xie2013sleep}. 

Investigation of the paravascular transport at macro-scale was proposed and tested in a rat's brain~\cite{iliff2013brain}. The procedure involved injecting
MRI contrast agent into the CSF and subsequently imaging the transportation of the MRI contrast agent at multiple time-points during a few hours after the injection. The MRI contrast agent worked as a CSF tracer, and was brain-wide in the rat after a few hours. The procedure was tested in humans for the first time 
in 2017~\cite{ringstad2017glymphatic} with acquisition of MRI images repeatedly during 48 hours after the injection and later quantified in a region-specific manner in 2018~\cite{ringstad2018brain} in individuals with dementia and controls. Overall, the MRI contrast agent transportation showed a centripetal pattern in all participants, but the MRI contrast agent was more protracted in individuals with dementia compared with controls~\cite{ringstad2018brain}. It was also noted that the CSF-tracer transport appeared faster than what can be expected from diffusion in simplified planar geometries. 

On this background, our purpose in this paper is to explore whether the CSF tracer distribution seen in~\cite{ringstad2018brain} can be explained by diffusion alone, as predicted by the seminal work of Sykov{\'a} and Nicholson~\cite{sykova2008diffusion}. 
We will investigate this hypothesis with finite element simulations of the diffusion process combined with a parameter identification procedure for the apparent diffusion coefficients (ADCs).
Thus, we aim to investigate whether we can assess ADC on long time-scales (hours or days), by fitting a diffusion model to the MRI data obtained at multiple time-points when the CSF tracer is propagating through the brain. This includes taking into account the complexity of the folding brain surface by constructing a patient-specific geometry. We propose that if the fit between images and our model is good and we identify ADC values that are in line with those predicted by diffusion tensor imaging (DTI) then enhanced solute transport of the MRI contrast agent in question can be ignored.  
Our approach for the parameter identification is to solve an optimization problem constrained by a diffusion equation with unknown coefficients, where the optimization targets the observed CSF tracer concentrations 
at the 10 available acquisitions during 24 hours after CSF tracer injection.


An outline of the paper is as follows: 
In Section~\ref{sec::methods}, we present the methodology of the paper. We start in Section~\ref{sec::method-mri} with a detailed description of the medical imaging methods relevant for this study. Section~\ref{sec::model} describes the mathematical models, and the computational methodology for this paper. 
In Section~\ref{sec::res}, we will present the results of the study, starting with the MRI analysis in Section~\ref{sec::res-MRI}. We continue in Section~\ref{sec::res-ass} with a synthetic test case, which involves finding robust regularization parameters with a uniform distributed noise added to the images.
The construction of the synthetic test case and the concentration estimation can be found in the Supplementary.
While in Section~\ref{sec::res-cont}, we present the computed ADC using the MRI images, and compare the values with ADC estimated with DTI. In Section~\ref{sec::res-pre}, we present the results of different method to decrease the boundary noise. The results will facilitate the general discussion in Section~\ref{sec::dis}.

\section{Methods}
\label{sec::methods}

\subsection{Simulation workflow}
\begin{figure}[h]
    \centering
    \includegraphics[width=0.575\textwidth]{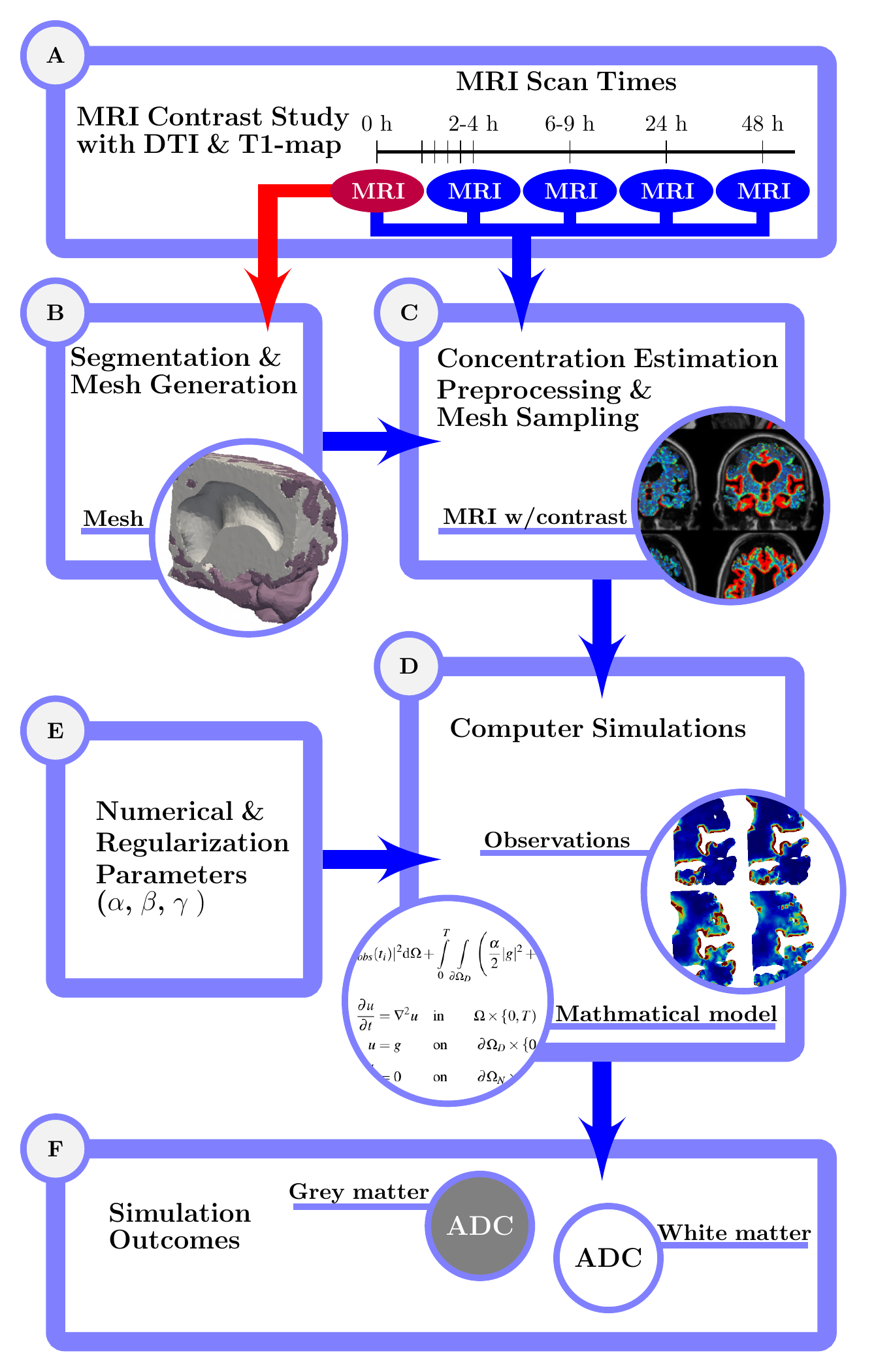}
    \caption{Schematic description of the simulation workflow, outlined through panels A to F.}
    \label{fig:flowchart}
\end{figure}

An overview of the simulation workflow for this paper is outlined in Fig.~\ref{fig:flowchart}. We obtained MRI data for one patient, which included MRI images with MRI contrast agent at different times (Box A in Fig.~\ref{fig:flowchart}). The first MRI image was segmented and used to construct a patient specific mesh (Box B in Fig.~\ref{fig:flowchart}). MRI images were used for estimating the CSF tracer concentration for the different times, and were subsequently sampled onto the patient specific mesh (Box C in Fig.~\ref{fig:flowchart}). The sampled concentrations at the different times were then used with the mathematical model (Box D in Fig.~\ref{fig:flowchart}). Values for numerical and regularization parameters were inputs for the computation (Box E in Fig.~\ref{fig:flowchart}). The simulations produced the optimal ADC for grey and white matter to explain the observations for different input parameters(Box F in Fig.~\ref{fig:flowchart}).

\subsection{Approvals and MRI Acquisition}
The approval for MRI observations was retrieved by the Regional Committee for Medical and Health Research Ethics (REK) of Health Region South-East, Norway (2015/96) and the Institutional Review Board of Oslo University Hospital (2015/1868) and the National Medicines Agency (15/04932-7). 
The study participants were included after written and oral informed consent.
The MRI images included 3D T1-weighted volume, sagittal 3D FLAIR, DTI and T1 map for the same patients.
All methods were performed in accordance with the relevant guidelines and regulations

The contrast observations were obtained using a 3 Tesla Philips Ingenia MRI scanner (Philips Medical Systems) with the same imaging protocol settings at all time points to acquire sagittal 3D T1-weighted volume scans. The imaging parameters were as follows: repetition time, “shortest” (typically 5.1 ms); echo time, “shortest” (typically 2.3 ms); flip angle, 8 degrees; field of view, 256  $\times$ 256 cm; and matrix, 256  $\times$ 256 pixels (reconstructed 512 $\times$ 512). We also obtained observation a sagittal 3D FLAIR volume sequence of the same patient, that was taken before the injection of contrast. The main imaging parameters were; repetition time = 4,800 ms; echo time 318 ms; inversion recovery time, 1,650 ms; field of view, 250 $\times$ 250 mm; and matrix, 250 $\times$ 250 pixels (reconstructed 512 $\times$ 512). The T1 map was obtained with a MOLLI5(3)3~\cite{taylor2016t1} sequence with the following imaging parameters; repetition time 2.3 ms; echo time 1.0 ms; flip angle 20; field of view 257 $\times$ 257; and matrix 240 $\times$ 240 pixels. The DTI acquisition was done with the parameters; repetition time 12171 ms; echo time 60.0 ms; flip angle 90;
field of view 240 $\times$ 240; and matrix 96 $\times$ 96 pixels.

\subsection{MRI Analysis}
\label{sec::method-mri}

\begin{figure}
\includegraphics[width=0.95\textwidth]{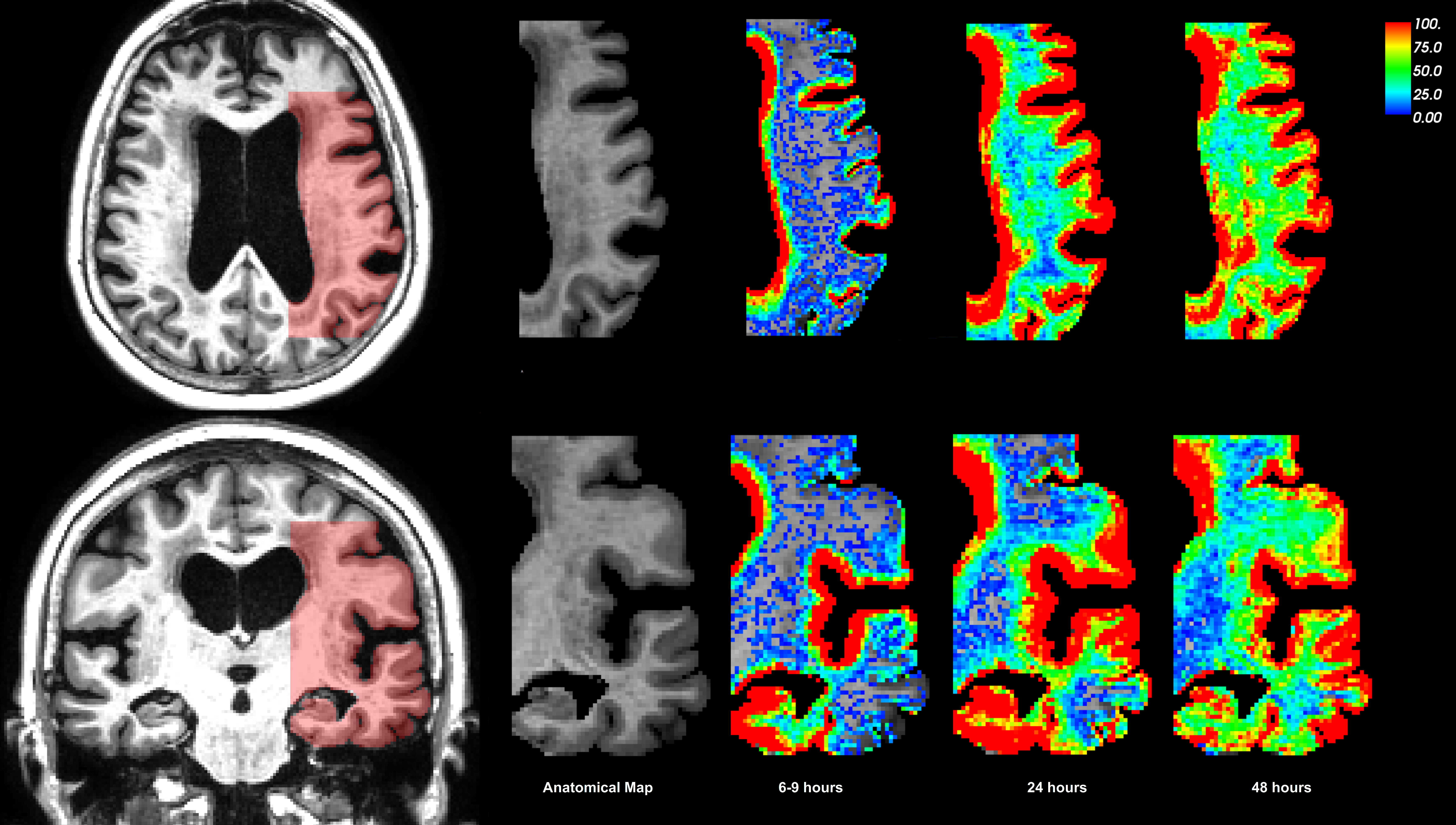} 
\caption{The image shows the percentage change in T1 signal unit ratios from baseline at different observation times in the slice (marked red in the left panel) used in the subsequent analysis. The color-bar was restricted to the range $(0,100)$. The upper row shows the axial slice, and the bottom row shows the coronal slice.}
\label{fig2} 
\end{figure}

A time sequence of T1-weighted MRI images showing the CSF- and brain enhancement in a patient diagnosed with NPH during 48 hours after intrathecal administration of Gadobutrol was obtained from a previous study~\cite{ringstad2018brain}. The software FreeSurfer \cite{Dale1999179, FischlLiuDale, spf2007, reuter:robreg10} was used to segment and align each of the observations, which made it possible to estimate voxelwise signal increase.
Figure~\ref{fig2} shows the distribution of MRI contrast agent in a selected region, as a percentage change in MRI signal unit ratios. The full data set used in this study (not all shown) consists of 10 MRI images, including a baseline MRI image taken before the contrast agent was injected. The MRI scans were obtained at different times distributed over 5 scans within 1-2 hours after injection, a single MRI scan at 2 hours, 6 hours, 24 hours and 48 hours. 
We segmented the baseline image with FreeSurfer and aligned the other images to the baseline. The exponential relation between the MRI signal values and the CSF tracer concentration, and the estimation of the concentration for each voxel is documented in the Supplementary~\ref{sec::singal-contrast}. The estimation of the concentration produced images similar to the MRI images, but the values have the unit millimolar (mM). Therefore, we will denote the concentration images as observations to distinguish the concentrations from the MRI image intensities.

\begin{figure}
\centering
\includegraphics[scale=0.2]{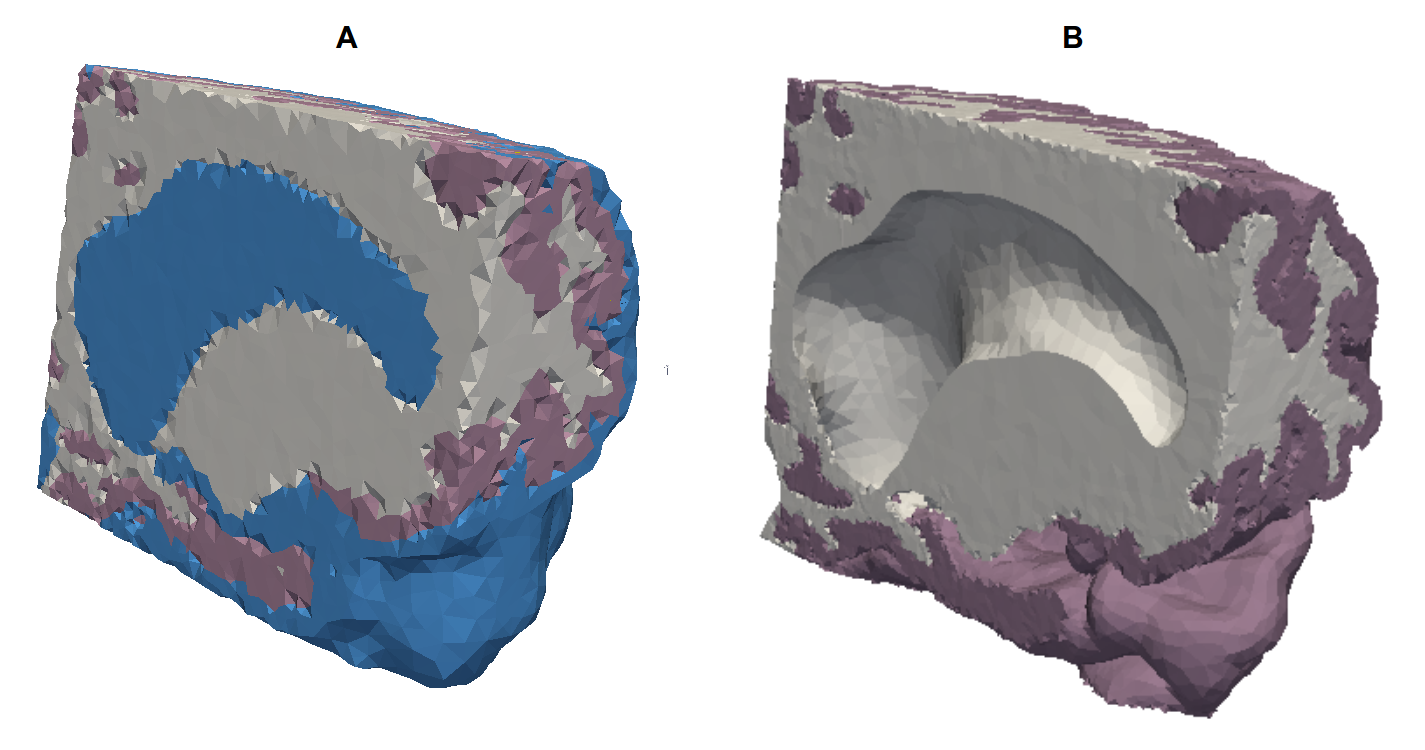} 
\caption{The leftmost image A) shows the mesh created from the baseline MRI image with three domains, while the rightmost image  B) shows the mesh created from the baseline MRI image with two domains. The blue domain corresponds to CSF domain, $\Omega_{CSF}$, the purple domain corresponds to grey matter, $\Omega_{GM}$, and the white domain corresponds to white matter, $\Omega_{WM}$. }
\label{Fig::Mesh}
\end{figure}

The segmentation process also produced polyhedral surfaces of the white and cortical grey matter that were used for mesh construction.
We used the Computational Geometry Algorithms Library (CGAL)~\cite{cgal:rty-m3-18b} to combine the surfaces
and construct the mesh with different subdomains. The computational requirement for the resulting mesh was significant, therefore two submeshes were also constructed, see Fig.~\ref{Fig::Mesh}.
The three domain mesh (white matter, grey matter and CSF compartment), shown in Fig.~\ref{Fig::Mesh} A, consists of 244,318 tetrahedral cells and 22,057 vertices, while the two domain mesh (without the CSF compartment), shown in Fig.~\ref{Fig::Mesh} B, consists of 183,138 tetrahedral cells and 42,514 vertices.


\subsection{Mathematical Model}
\label{sec::model}
The macroscopic extracellular diffusion in the brain can be considered a hindered diffusion with an apparent diffusion coefficient (ADC) depending on the structure of the extracellular space\cite{sykova2008diffusion}. 
The relation between the apparent and free diffusion coefficients is defined as 
\begin{equation}
 \lambda =  \sqrt {D/D_{ADC}}
\label{tortuosity}
\end{equation}
with $\lambda$ denoting the tortuosity of the extracellular space.
In order to estimate the ADC involved in the contrast transportation, shown in Fig.~\ref{fig2}, we assume that the process can be modeled by a diffusion equation. 
Then, we constructed an optimization problem with the aim of minimize the difference between the observed and the modeled contrast distribution by optimizing the boundary conditions and the apparent diffusion coefficient. Thus enhanced transportation because of effects such as dispersion would result in an ADC larger than that predicted by DTI. 
The optimization problem was defined as 
\begin{equation}
\min_{D,g} \quad \sum\limits_{i=1}\sp{n} \int\limits_{\Omega} |u(t_i) - u_{obs}(t_i)|\sp{2} \mathrm{d}\Omega + \int\limits_{0}\sp{T} \int\limits_{\partial \Omega_D} \left( \frac{\alpha}{2} | g |\sp{2} + \frac{\beta}{2} \left| \frac{\partial g}{\partial t} \right|\sp{2} +  \frac{\gamma}{2}| \nabla g |\sp{2} \right) \mathrm{d}\Omega \mathrm{d}t  
\label{EQ::objf}
\end{equation}
subject to   
\begin{equation}
\begin{aligned}
\frac{\partial u}{\partial t} &=  D \nabla\sp{2} u && \text{in} \qquad \Omega \times \left\lbrace 0 , T \right)  \\
u&=g && \text{on} \qquad \partial\Omega_D  \times \left\lbrace 0 , T \right) \\
\frac{\partial u}{\partial n}&=0 && \text{on} \qquad \partial\Omega_N  \times \left\lbrace 0 , T \right) 
\end{aligned}
\label{Eq::PDE}
\end{equation}
Here, $u \left[ \mathrm{mM} \right] $ is the simulated, time-varying CSF tracer distribution, $D \left[ \mathrm{mm^2/s} \right]$ is the ADC, $g \left[ \mathrm{mM} \right] $ is the boundary condition, $\Omega$ is the domain, and $T \left[ \mathrm{h} \right] $ is the final simulation time. We assume that the domain $\Omega$ consists of three sub domains, each with a different ADC. We denote the CSF (subarachnoid and lateral ventricle) domain as $\Omega_{CSF}$, the grey matter as $\Omega_{GM}$ and the white matter as $\Omega_{WM}$. 
The ADC was assumed to be constant within the CSF, grey and 
white matter but each region may have different values. The $\alpha$, $\beta$ and $\gamma$ parameters are non-negative regularization parameters and $u_{obs} \left[ \mathrm{mM} \right] $ are the concentration distribution at time-points $t_i  \left[ \mathrm{h} \right]$. 
Spacial regularization parameter $\alpha$ enforces smoothness on the boundary by minimizing the concentration, i.e. high value of $\alpha$ will give less concentration in the optimal solution. Temporal regularization parameter $\beta$ enforces smoothness in time on the boundary, i.e. high value of $\beta$ will give a smoother concentration curve in time. Gradient regularization parameter $\gamma$ enforces continuity between adjacent concentrations at  the boundary, i.e. high values of $\gamma$ will give smoother concentration values at the boundary.

\subsubsection{Boundary conditions.}

For the three domain geometry (with grey and white matter, and CSF compartment), the Dirichlet boundary condition $\Omega_D$ was only applied on the outward facing boundary of the CSF domain, $\partial \Omega_1$. Homogeneous Neumann conditions $\Omega_N$ were applied on the remaining boundaries.

\begin{figure}

\centering
\includegraphics[width=0.4\textwidth]{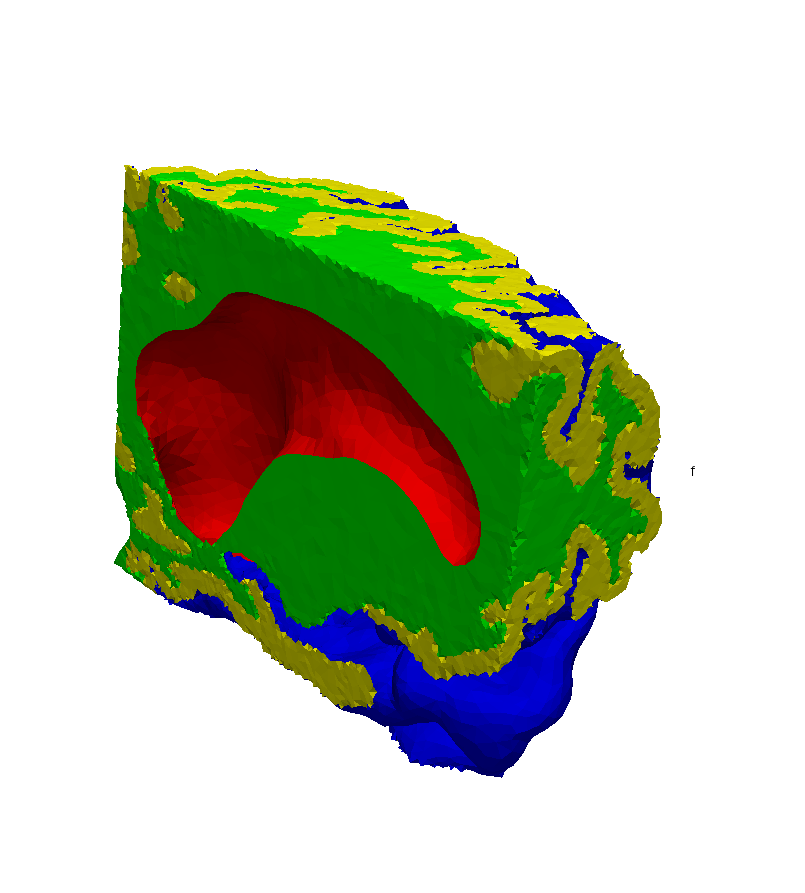} 
\caption{The images shows the different boundary values on the mesh, with each color representing a unique boundary. }
\label{markedboundaryvalues}
\end{figure}

The implementation of the gradient regularization $\gamma$ for the case containing only grey and white matter required that the outward facing boundary was decomposed in different regions to avoid the boundary values being continuous at the interface between CSF, grey and white matter. We decomposed the boundary as seen in Fig.~\ref{markedboundaryvalues}, with the red and blue boundary adjacent to the CSF. We defined the red boundary $\partial \Omega_r$ and the blue boundary $\partial \Omega_b$ as Dirichlet boundaries  $\Omega_D$, while the green and yellow were Neumann boundaries  $\Omega_N$. The regularization parameter $\gamma$ was subsequently set to be non-zero in~\eqref{EQ::objf}. We initially tested gradient regularization with the same parameter $\gamma$ on both boundaries, but the different distribution of tracers on the boundaries made it difficult to find an adequate value for $\gamma$. This may be attributed to the fact that the concentration in the lateral ventricles were more uniform than that in the SAS, so we defined $\gamma$ as
\begin{equation}
\gamma = \begin{cases}
0.01\tilde{\gamma}  \in &\partial \Omega_{b} \\ 
\tilde{\gamma} \qquad \in  &\partial \Omega_{r}
\end{cases}  
\end{equation} 
with $\tilde{\gamma}$ as the referenced, i.e. mentioned in the text, regularization parameter. 


\subsubsection{Synthetic test case with a manufactured solution.}

In order to assess the robustness and accuracy of the methodology of ADC estimation via PDE constrained optimization we constructed a synthetic test case with a known, manufactured solution. The setup for the numerical tests can be found in the Supplementary~\ref{sup:manu}.  In the case of three domains we varied $\alpha \in (10^{-6}, 10^{-2})$ and $\beta\in(10^{-6}, 10^2)$. In the case of two domains $\alpha\in(10^{-6}, 1)$, $\beta\in(10^{-4}, 10^2)$ and $\gamma\in(10^{-4}, 1.0) $. 

We tested the noise susceptibility with a uniform distributed of noise. This was done by adding noise in the range of $(-n_{amp}, n_{amp})$ to the observation after loading, i.e. each vertex. We tested the noise with $n_{amp}\in(0.03,0.15,0.30,0.475)$, compared to the maximum values for the manufactured solution being in the range (0.3,1.3) for all times. We  used the observation with the added noise to compute and confirm the expected SNR, since there was a finite number of vertices. The variation in the number of observations was tested with 5,10 and 20 evenly spaced observations over the course of 24 hours combined with 10 times steps, 20 time steps and 40 time steps. 

\subsection{CSF tracer distribution reconstruction from MRI data}

\label{sec::met::contrast-dist}
The MRI data consisted of scans at times $t_i$ that were distributed in the 48 hours time frame, with the first observation with tracer occurring 1-2 hours after the injection. This was followed by 4 observations within the first hour, and we observed no visible change in the tracer distribution for these observations. Therefore, we used observation times listed in Tab.~\ref{model-params-overview} for the computation of the ADC values. 


\subsubsection{Comparison with data obtained from DTI analysis}
\label{sec::met::compare}
We compared the ADC computed by solving \eqref{EQ::objf}-\eqref{Eq::PDE} with ADC values obtained from the DTI image of the same patient. We remark that a direct comparison with DTI is not possible because DTI measures the ADC of water. Therefore, we used the DTI to estimate the tortuosity in grey and white matter, which together with the free diffusion coefficient of Gadobutrol and~\eqref{tortuosity} can be used to approximate ADC for Gadobutrol, details and references are found in the Supplementary \ref{sec:supp:dti}.


\subsubsection{Preprocessing of concentrations of Gadobutrol at ventricular and subarachnoid surfaces}
High frequency concentration changes was observed at the boundaries of our mesh, see top row of column A in Fig.~\ref{SAS} and Fig.~\ref{VENT}, which can be interpreted as sampling errors. Such errors may be caused by noise in the MRI data, errors in the segmentation, the segmented polyhedral surfaces which typically cut voxels, miss-alignment between different observations, and by the inaccuracy of sampling discontinuous voxel data. Therefore, we investigated two approaches to reduce high frequency components at mesh boundaries:
\begin{itemize}
\item A projection of the segmented CSF Gadobutrol concentrations directly at the ventricle and subarachnoid surfaces (CP), 
\item A Gaussian smoothing (GS) procedure.
\end{itemize} 
In detail, the CP method was implemented by finding the voxel corresponding to each boundary vertex on $\Omega_b$ and $\Omega_r$ using an affine transformation matrix in FreeSurfer. Then, for each voxel, we computed the average of all surrounding voxels with CSF segmentation mark in a 7$\times$7$\times$7 matrix. The average values were then used at the corresponding boundary vertex in the computations. In order to estimate the CSF concentration, we assumed that the T1 value corresponding to CSF was 3000 ms. Finally, for the GS method we used the Gaussian smoothing function found in the python-module scipy~\cite{jones2001scipy}, and applied the smoothing to all voxel in the  observations estimated from MRI images. The standard deviation of the Gaussian distribution was set to 1.5 mm$\sp{2}$ compared to MRI voxel length of 1.0 mm. In addition to the CP and GS methods, we also use the raw data without any
preprocessing. This method is referred to as RAW. 

We performed the simulations with the regularization parameters $\alpha\in(10^{-6},10^{-4}), \beta\in(1.0,10) , \gamma\in(0.0,0.01,1.0)$. Additionally, we ran simulation with $\beta=100$ and $\gamma=100$ to obtain simulations that clearly showed differences in the concentration distributions on the boundary.

\begin{table}[hhhh]\centering
\begin{tabular}{|ccc|}
\hline
Description & Model variable  & Value used in the simulation\\
\hline
Number of timesteps &  k	 & 24, 48	\\
\hline
Observational time points  & $\tau$ &  $\lbrace 0.0 ,2.1 , 6.1 , 24.1, 47.8 \rbrace$	 \\
\hline
Timestep & $dt$	 	   &	  1, 2 $\mathrm{hours}$	\\ 
\hline
Spacial Regularization & $\alpha$	   &	  $\lbrace 10\sp{-6}, 10^{-4}\rbrace$\\ 
\hline
Time Regularization   & $\beta$	   &	 $\lbrace 1.0 , 10.0\rbrace$	\\ 
\hline
Gradient Regularization   & $\beta$	   &	 $\lbrace 0.0, 0.01, 1.0\rbrace$	\\ 
\hline
\end{tabular}
\caption{The table shows an overview of the parameters used for boundary constrained optimization with MRI data. }
\label{model-params-overview}
\end{table}

\section{Results}
\label{sec::res}
\subsection{MRI Analysis}
\label{sec::res-MRI}
Statistical analysis was done on the observations to determine the signal to noise ratio (SNR) for each observation. The first 5 observations had SNR in the range 0.6-0.8 in grey matter and 0.2-1.3 in white matter. For the remaining observations the SNR ranged 1.4-1.9 in grey matter and 1.3-2.6 in white matter.

\subsection{Assessment of accuracy and robustness on a synthetic test case}
\label{sec::res-ass}



For the first geometry, involving CSF, grey and white matter, we ran a series of 448 tests with different regularization parameters $\alpha$ and $\beta$ using the manufactured solutions with observations every 2.4 hours over the course of 24 hours. It was found that for $\alpha \in (10^{-6}, 10^{-2} ) $ and $\beta\in(10^{-6}, 1.0)$, the error in
the ADC for CSF, grey and white matter was less than 5\%. For combinations with $\beta=10^2$, the smallest error was 63.6\% in the CSF, 10.8\% in grey matter and 3.5\% in white matter. While for $\alpha=1.0$, the smallest error was 160.2\% in the CSF, 23.2\% in grey matter and 8.9\% in white matter. Second, the robustness of the parameter identification process with respect to noise in the data was investigated. The noise values were randomly obtained from a uniform distribution in the range (-0.3,0.3), where 0.3 equaled the maximum initial value of the manufactured solution
and was 23\% of the manufactured solution at its max. Again, for $\alpha \in (10^{-6}, 10^{-2})$ and $\beta\in(10^{-6}, 1.0)$ the error in the ADC in CSF was less than 23\% and less than 9.7\% in grey and white matter. 

For the second geometry, including only grey and white matter, we ran a selection of 186 tests to ensure that the results were consistent. The boundary conditions were applied to the boundaries of the SAS and lateral ventricle. In the range $\alpha\in(10^{-6}, 10^{-2}) $ and $\beta\in(10^{-4}, 1.0) $, the error was less than 4.2\% for the ADC in both grey and white matter. 

The noise susceptibility was tested with the addition of the regularization parameter $\gamma\in(0.0,10^{-2},1.0)$, which enforces smoothness at the boundary. We computed the error on the SAS boundary and on the lateral ventricle boundary, and it was observed that $\gamma$ did not contribute to a lower error in the ADC, but decreased the boundary error with a few percentage on average. In the range $\alpha\in(10^{-6}, 10^{-4}) $ and $\beta\in(10^{-4},1.0)$, we had maximum error of 7.0\% in the grey matter and 3.1\% in the white matter with noise randomly obtained from the uniform distribution range of (-0.3,0.3). The computed the SNR of 0.74 for first time step, and increased to the maximum of 6.3.

The synthetic test case revealed that the second geometry, involving only the estimation of ADC in grey and white matter was the most accurate method. Hence, in the following, only the second geometry was used.



\subsection{CSF tracer distribution recondstrction from MRI data}

\label{sec::res-cont}
\begin{figure}
\centering
\includegraphics[width=0.8\textwidth]{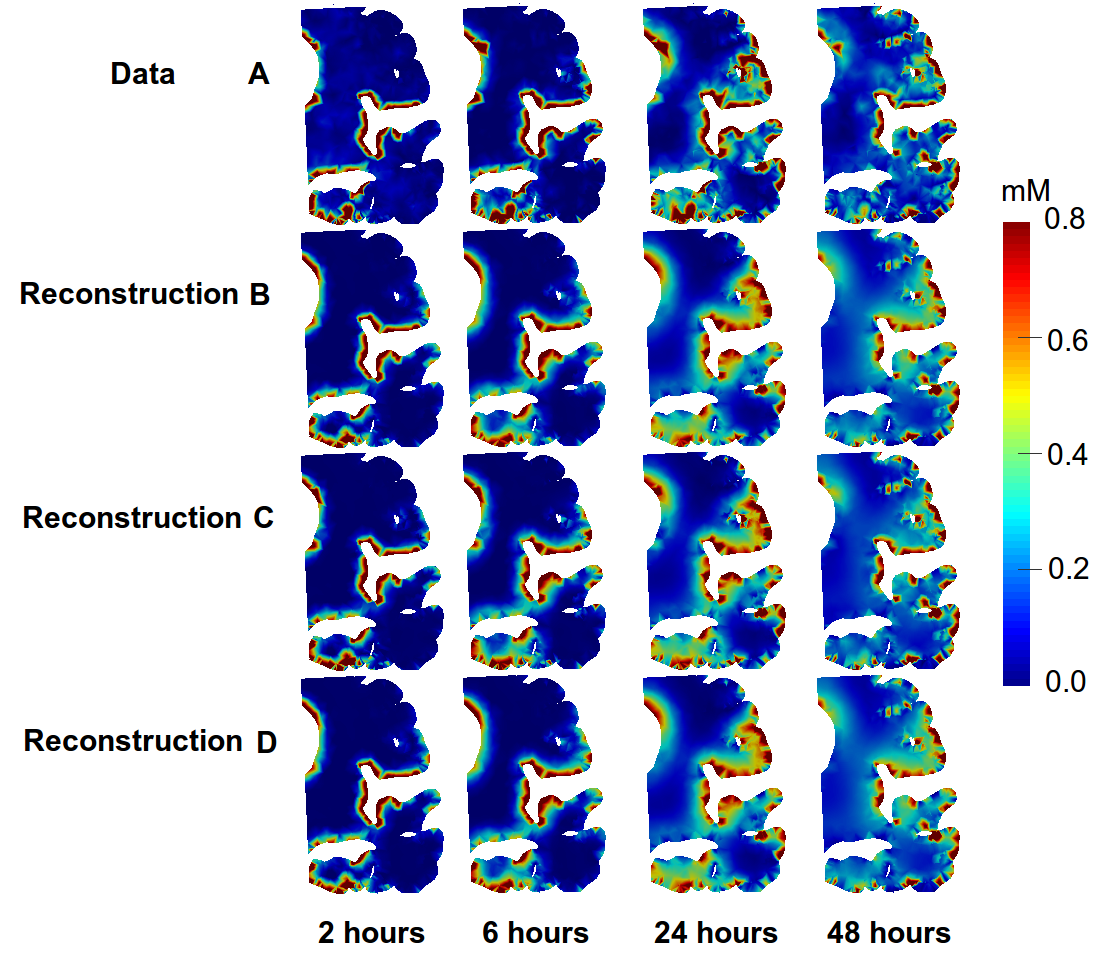}
\caption{
The image displays observations from MRI images and computational reconstruction of the same observations at observational time-points, i.e. 2 hours, 6 hours, 24 hours and 48 hours after the administration of CSF tracer. Row A) shows the observations estimated from MRI images. Row B) shows the reconstructed observations with $\alpha=0.0001,\beta=10.0,\gamma=1.0$ and 48 time steps. Row C) shows reconstructed observations with $\alpha=10^{-6},\beta=1.0,\gamma=0.01$ and 48 time steps. Row D) shows reconstructed observations with $\alpha=10^{-6},\beta=10.0,\gamma=1.0$ and 48 time steps. The color-bar was restricted to the range (0.0-0.8). 
 }
\label{Fig::realdata}
\end{figure}

In the synthetic case we managed to reproduce the ADC within 7\% error for a wide range of parameters, $\alpha\in(10^{-6}, 10^{-2}) $, $\beta\in(10^{-4}, 1.0) $, and $\gamma\in(0.0,10^{-2},1.0)$. The same parameters were then used to compute the optimal ADC for the CSF tracer estimated from MRI images. A few of the reconstructions are shown in Fig.~\ref{Fig::realdata}. It is clear that our model assumption of an underlying diffusion equation is adequate, that is; from a visual point of view the observed data in Fig~\ref{Fig::realdata}, row A, is reconstructed accurately for the various 
different regularization parameters Fig~\ref{Fig::realdata}, row B-D. 
We computed the average ADC over the range of regularization parameters to be 0.57$\pm$0.05 mm$^2$/h in grey matter and 0.72$\pm$0.02 mm$^2$/h in the white matter, which respectively corresponds to $1.6\pm 0.2\times 10\sp{-4} \mathrm{mm^{2}/s}$ and $2.0\pm 0.1 \times 10^{-4} \mathrm{mm^{2}/s}$, shown in Fig.~\ref{diffcoeff}.



\subsubsection{Comparison with data obtained from DTI analysis}
\label{sec::res-comp}
The median ADC for water in the DTI was estimated to be
$1.0\pm0.4 \times 10\sp{-3} \mathrm{mm^{2}/s}$ in grey matter and
$0.9\pm0.3 \times 10\sp{-3} \mathrm{mm^{2}/s}$ in white matter 
resulting in a tortuosity of 1.73 and 1.85 based on~\eqref{tortuosity}. The free diffusion coefficient for Gadobutrol was approximated to be $3.8\times 10^{-4} \mathrm{mm^{2}/s}$ in the Supplementary \ref{sec:supp:dti}. This gives an estimate median Gadobutrol ADC to be 1.3$\pm0.5\times 10^{-4}\mathrm{mm^{2}/s}$ in the grey and 1.1$\pm0.4\times 10^{-4}\mathrm{mm^{2}/s}$ in the white matter. This estimation assumed that the tortuosity was independent for molecules with mass lower than 1kDa. For Gadobutrol, this gives an  average of 23$\pm$15\% larger ADC in grey matter and 82$\pm$10\% larger ADC in white matter compared to the associated ADC values based on DTI. 

\subsubsection{Preprocessing of concentrations of Gadobutrol at ventricular and subarachnoid surfaces }
\label{sec::res-pre}

\begin{figure}
\centering
\includegraphics[scale=0.35]{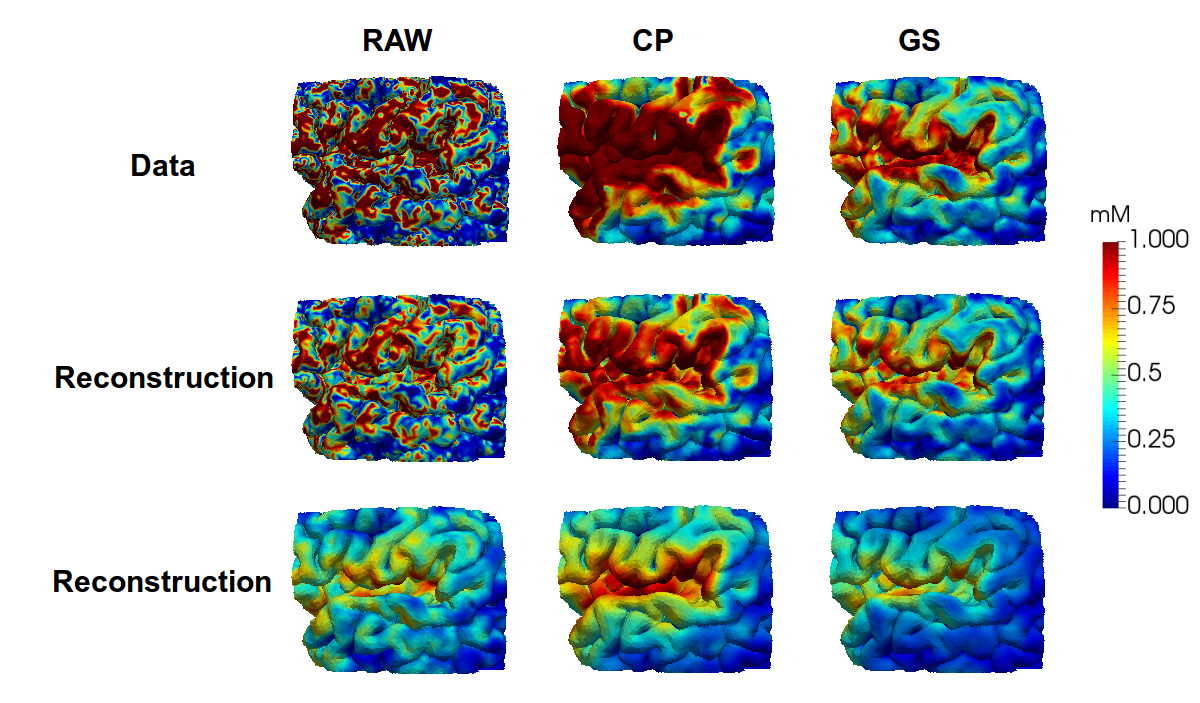}
\caption{ The images shows the SAS boundary after 12 hours. The upper row shows the observations with following preprocessing left to right: Raw observations, projection of CSF value onto the boundary, Gaussian smoothing. The middle row shows the corresponding states with the regularization parameters $\alpha , \beta, \gamma = (10^{-6},1.0 ,1.0)$. The bottom row shows the corresponding states with the regularization parameters $\alpha , \beta, \gamma = (10^{-6},100.0 ,100.)$. }
\label{SAS}
\end{figure}

\begin{figure}
\centering
\includegraphics[scale=0.35]{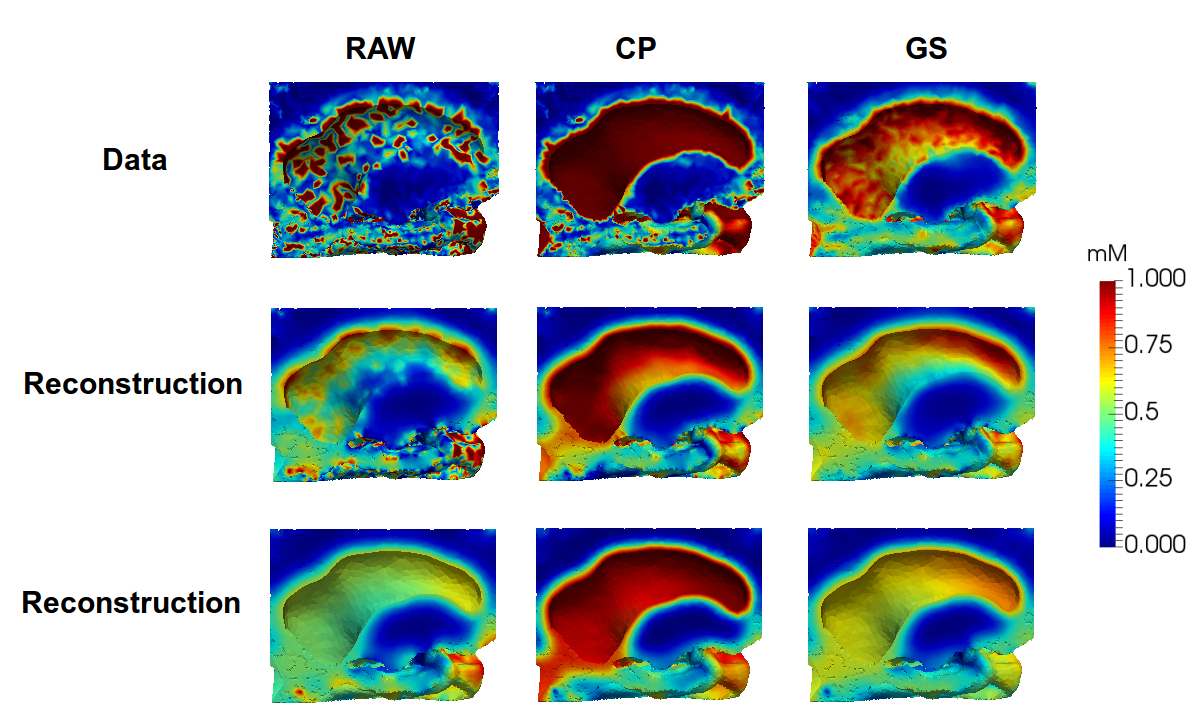}
\caption{ The images shows the ventricular wall after 12 hours. The upper row shows the observations with following preprocessing left to right: Raw observations, projection of CSF value onto the boundary, Gaussian smoothing. The middle row shows the corresponding states with the regularization parameters $\alpha , \beta, \gamma = (10^{-6},1.0 ,1.0)$. The bottom row shows the corresponding states with the regularization parameters $\alpha, \beta, \gamma = (10^{-6},100.0 , 100) $ .}
\label{VENT}
\end{figure}

\begin{figure}
\centering
\includegraphics[scale=0.3]{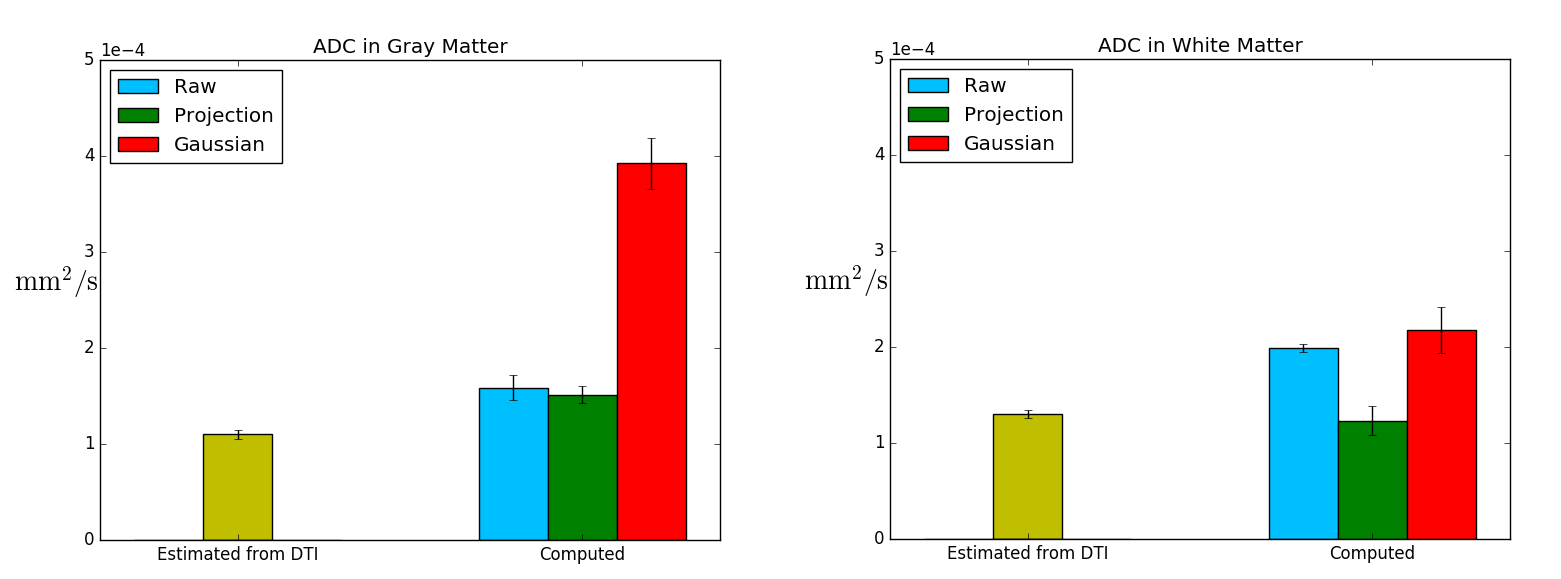}
\caption{
The images shows the average ADC estimated with DTI and the computed ADC in grey and white matter for different sampling methods and regularization parameters $\alpha\in(10^{-6},10^{-4}), \beta\in(1.0,10.0), \gamma\in(0.0,0.01,1.0)$ and number of time steps $k\in(24,48)$. The error bars show the standard deviation.}
\label{diffcoeff}
\end{figure}

The three different reconstructions with RAW, CP and GS methods are illustrated in Fig.~\ref{SAS} and Fig.~\ref{VENT} for some different regularization parameters.
As can be seen in RAW column in Fig.~\ref{SAS} and Fig.~\ref{VENT} the boundary gradient regularization $\gamma$ caused the concentration to be more uniform on the boundary and in particular for the high value $\gamma=100$. The CP method results in similar values for most of the ventricular boundary. The resulting ADC Gadobutrol values for RAW, CP and GS are shown in Fig.~\ref{diffcoeff} together with the Gadobutrol ADC value estimated from  DTI. We computed the percentage difference for each computed ADC and compared with ADC estimated from DTI. For the CP method, the average difference was 37$\pm$8\% in grey matter, -5$\pm$11\% in white matter and for the GS method, the average difference was 250$\pm$22\% in grey matter, 68$\pm$18\% in white matter.

\section{Discussion}
\label{sec::dis}
The glymphatic system proposes that the paravascular network facilitates brain-wide transportation by bulk flow. Several mathematical
modeling studies have been performed at the micro-level, but to our knowledge, our current study
is the first study to investigate this process based on 
human imaging data at the macro-level. Our investigations are based on the distribution of Gadobutrol in CSF and brain tissue up to 48 hours after intrathecal injection. Our results here are supportive of the glymphatic system in the sense that the ADC estimated from the Gadobutrol distribution 
is higher than corresponding numbers obtained from DTI. Hence, our study suggests that at longer time scales, the slow glymphatic system has an impact.  


We assumed that the white matter was isotropic, however, it is well-known that the white matter is anisotropic. The main reason for our choice was that introducing anisotropy would result in a model with many more free variables, requiring substantial parameter tuning, and consequently result in less predictability. Furthermore, it
can be seen in Fig. \ref{Fig::realdata} that the reconstruction is quite good. Demanding that the transport in the CSF was governed
by a diffusion process yielded a much larger error in the synthetic test case and this assumption should probably be avoided.

In our approach we have employed the FreeSurfer toolkit to segment and register grey and white matter surfaces. FreeSurfer provides segmentations with sub-voxel accuracy, which 
means that at the boundaries towards the ventricles and subarachnoid space, the surface boundary is usually in the interior of a voxel rather than at the voxel 
boundaries. The consequence is that the raw data appear to have noisy image intensities at the boundaries which cut voxels. For this reason, we investigated two different approaches
to interpret the data at the boundaries, in addition to using the raw data itself. 
We observed in Fig.~\ref{diffcoeff} that the GS method increased the computed ADC with approximate 250$\pm 22\%$ in grey matter, while the CP method corresponded best with -5$\pm$11\% difference to the ADC estimated with DTI in white matter. However, it should be noted that the CP method imposes concentration values in the CSF onto the tissue boundaries, not accounting for the transverse propagation of MRI contrast agent through a membrane, like the pia mater or the endothelial layer of the ventricles. Hence, potentially, an unnaturally large concentration gradient at the boundary that is caused by a partial barrier rather than the white matter itself may have resulted in the low ADC seen when using this method. 
Which of the methods that best depicted the actual boundary concentration is unknown
and would need to be determined by phantom studies. 
Another limitation is that we have not yet been able to assess the diffusivity of Gadobutrol and have  have relied on literature values for  the free diffusion coefficient of a similar mass molecule as a substitute for Gadobutrol. However, it should be noted that even though the concentration values appear noisy at the boundaries, Fig.\ref{SAS}, the interior reconstructions appear accurate, Fig. \ref{Fig::realdata}.

Concerning the imaging, a current limitation is that standard T1 weighted MRI images have higher resolution and
SNR than corresponding T1-maps. In detail, the intrathecal contrast enhanced T1 weighted volume scans had 1.0 mm slice thickness, while the T1-map slice thickness was 4.0 mm. This means that the calculation of the concentration at boundaries can suffer from mismatch of tissue and CSF. Furthermore, the T1-map sequence is designed to estimate the T1 times in tissue, and therefore does not give accurate values for the CSF.  The average T1 relaxation time for the CSF in left lateral ventricle was $181\pm349 \mathrm{ms}$, compared to the value 1000-5500 ms that can be found in the literature~\cite{condon1987mr} and we used the literature values to compute the concentration in the CSF. We also estimated the T1 relaxation time for grey and white matter to be 1200$\pm$271 and 819$\pm$180 ms, which compares better with the literature values ranging 1470-1800 ms in grey and 1084-1110 ms in white matter~\cite{stanisz2005t1}.

The ADC of water in this study was measured with DTI to be $1.0\pm0.4 \times 10\sp{-3} \mathrm{mm^{2}/s}$ in grey matter and $0.9\pm0.3 \times 10\sp{-3} \mathrm{mm^{2}/s}$ in white matter. 
It has been reported that ADC of extracellular water in young and healthy subjects 
was $0.78-1.09 \times 10^{-3}$ in the cortical grey matter and
$0.7-0.9\times 10^{-3} \mathrm{mm\sp{2}/s}$ in white matter~\cite{helenius2002diffusion}. Although our values match with the DTI values of this study, the white matter seem to be on the upper threshold. This can be explained by the fact that subjects with dementia typically have higher ADC values, less anisotropy and greater variation in the white matter\cite{goujon2018can}. Thus, the DTI values found in this study seems to correspond with the literature values. 
Based on the raw data, our estimated Gadobutrol ADC was $23\pm 10$\% and $82\pm 4$\% larger in the grey and white matter, respectively, than what
the DTI data corresponds to for Gadobutrol.

Previous mathematical modeling studies at the micro-level~\cite{holter2017interstitial, smith2017glymphatic} suggest that diffusion dominates in the interstitium. However, diffusion depends on molecular size as described by the Stokes-Einstein equations and large molecules
are transported slower than small molecules~\cite{valencia2011understanding}.
Convective flow of solutes, on the other hand, is independent of the molecular size. Furthermore, transport has been reported to be independent of molecular size~\cite{cserr1981efflux}, a fact that suggests convective transport. In fact, convective velocities of $0.8-4\times 10^{-3}\mathrm{mm/s}$ has been demonstrated or estimated~\cite{cserr1991extracellular,kress2014impairment,ray2019analysis,xie2013sleep}, indicating that the solute transport would be dominated by convection for large molecules, whereas similar to diffusion for smaller molecules, such as water. Gadobutrol is in this context a molecule of moderate size, i.e., 604 Da, and hence 
not ideally suited for the study of bulk flow for larger molecules such as A$\beta$  (4.5 kDa) or CSF-$\tau$ (45 kDa). In fact, Gadobutrol is predicted to have 
a Pechlet number less than one\cite{holter2017interstitial,ray2019analysis} which from these micro-level studies would imply that the distribution of
Gadobutrol is governed mainly by diffusion.

Studies~\cite{asgari2016glymphatic, brynjfm, Diem} have found that dispersion in the paravascular spaces adds less than a factor two to diffusion for solute transportation. However, all these studies were done with modeling that was on the micro-scale over shorter time periods. To the authors' knowledge, the only other study~\cite{croci2019uncertainty} that has considered macroscopic modeling on the time-scales of hours and days, where uncertainties representing both variations in ADC and paravascular velocities where modelled with extensive testing using Monte Carlo methods. They found that, in particular, the CSF tracer distribution within the deep white matter found in~\cite{ringstad2018brain} could not be explained by 
diffusion alone. Recent high-resolution 
MRI imaging of paravascular spaces in a rat's brain points towards significant contributions from white matter paravascular spaces connected to the ventricles\cite{magdoom2019mri}. 
It should also be mentioned that the permeability of the white matter is several orders of magnitude higher than the grey matter~\cite{dai2016voxelized,stoverud2012modeling}
and may as such be more  susceptible to bulk flow than grey matter. 


In conclusion, we computed that the ADC in grey and white matter and found that in both cases the ADC was somewhat larger than estimates based on DTI alone. Thus, indicating that there is potential for enhanced solute transportation in the brain over a longer time period. There are, however, a number of uncertainties that needs to be taken into account, for instance the resolution of DTI and T1 mapping.

\subsection*{Author contributions statement}
L.M.V, P.K.E, G.R., KAM conceived the experiments.
L.M.V., S.K.M, S.F. implmented the simulators.
L.M.V conducted the experiments and made the figures.
All authors discussed and analyzed the results. 
L.M.V, K.A.M wrote the first draft. 
All authors revised the manuscript
and approved the final manuscript. 

\subsection*{Additional information}
\textbf{Competing interests} The authors declare no competing interest. 

\subsection*{Acknowledgement}
The simulations
were run on the Abel Cluster (Project NN9279 K), owned by the University of
Oslo and the Norwegian Metacenter for High-Performance Computing and
operated by the Department for Research Computing at University Center
for Information Technology, the University of Oslo Technical Department,
www.hpc.uio.no/. 

\subsection*{Corresponding Author}
The corresponding author's e-mail is kent-and@math.uio.no. 

\subsection*{Data Availability}
The datasets analyzed in the current study are available from the corresponding author upon request.

\bibliography{references}

\clearpage
\section{Supplementary}

\subsection{Constructing the synthetic solutions to assess accuracy and robustness.}
\label{sup:manu}

Below we will discuss the parameter identification and its sensitivity with respect to the regularization parameters, noise, number of observations and time-resolution of the forward model. 
The manufactured observations used in \eqref{EQ::objf} were obtained by forward computation of \eqref{Eq::PDE} with the Dirichlet boundary condition defined as 
\begin{equation}
g(t) = 0.3 +0.167t - 0.007t^{2} \qquad \text{ for } 0 \leq t \leq 24.
\label{EQ::DIRI}
\end{equation}
The initial condition was set to 0 everywhere, the time step was $dt = 0.24$, and the diffusion coefficients were selected to be 
\begin{equation}
D_{\Omega_1} = 1000.0, \quad D_{\Omega_2} = 4.0, \quad D_{\Omega_3} = 8.0 
\end{equation}  

\begin{figure}[b]
\centering
\includegraphics[scale=0.4]{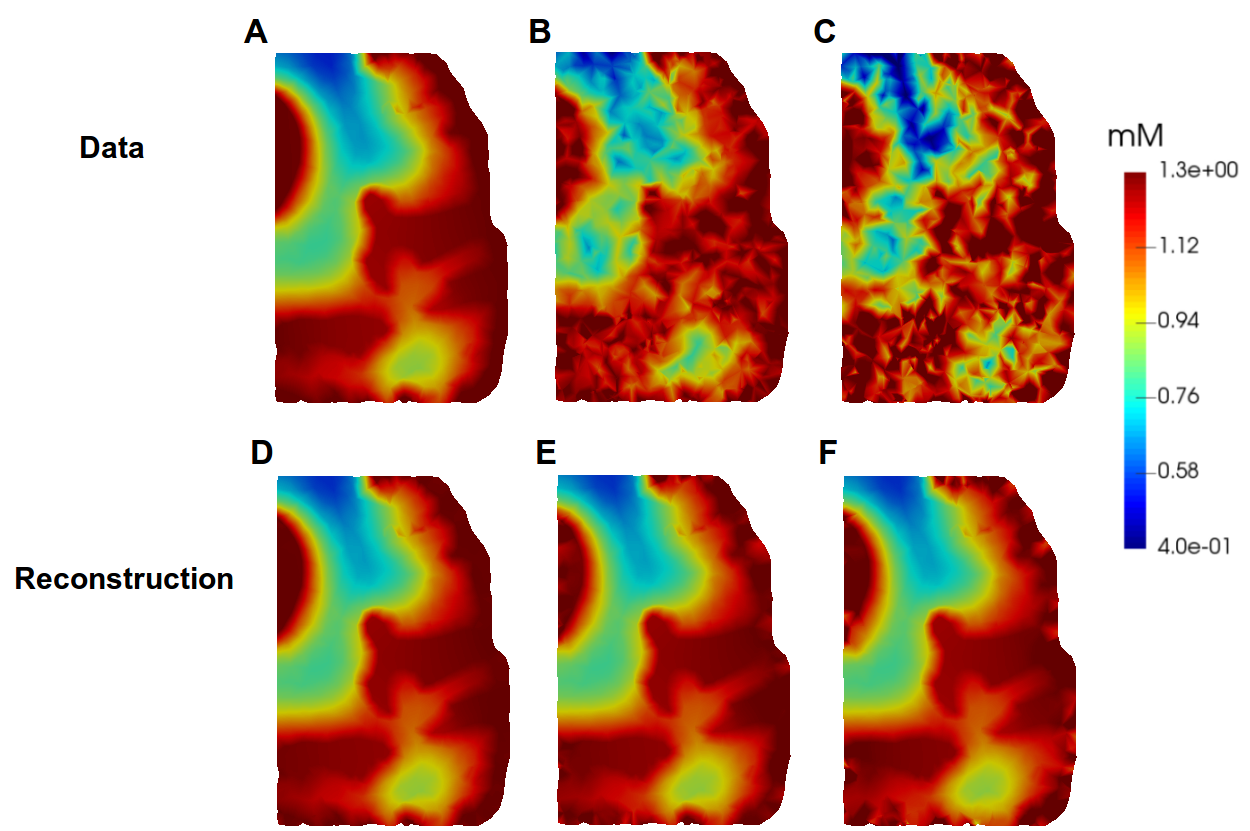}  
\caption{The upper row shows the manufactured observation, A) Shows the manufactured observation at time-point 24 with no noise added. B) Shows the manufactured observation after 12 hours with noise range of (-0.15,0.15). C) Shows the manufactures observation after 12 hours with noise range of (-0.3,0.3). The lower row shows the results with optimized parameter obtained with $\alpha=0.0001$, $\beta=1.0$ and $k=20$. D) Shows the resulting state given the observation in A. E)  Shows the resulting state given the observation in B .F) Shows the resulting state given the observation in C. }
\label{12hourswithnoise}
\end{figure}

\subsection{Additional results: Assessment of accuracy and robustness on a synthetic test case}
We will in this section present additional results that was not presented in the main text, Sec. \ref{sec::res-ass}. The convergence was monitored for a wide range of parameters, some shown in Figure~\ref{convergence}, and we saw the similar convergence for ADC and boundary conditions for $\alpha\in(10^{-6}, 10^{-4})$ and $\beta\in(10^{-4}, 1)$. 
Furthermore, the method was found robust with respect to the noise as illustrated in the Figure~\ref{12hourswithnoise}.

The variation in the number of observations was tested with 5,10 and 20 evenly spaced observations over the course of 24 hours. This was done together with 10 times steps, 20 time steps and 40 time steps. Illustrative results are shown in Figure~\ref{boundarycontrol}, which shows that a high number of time steps compared to observations causes oscillations (e.g. purpule line) at the boundary for $\alpha\gg\beta$. This is counteracted by selecting high values of the temporal regularization $\beta$ relatively to $\alpha$ hereby enforcing smoothness (e.g. blue line). For example, in the geometry with grey and white matter, 
it was observed that the ADC error increased by a factor $\approx 4$ for $\alpha\sim\beta$  relative to the case where $\beta\gg\alpha$ for  the case with 20 time steps and 10 observations.

\subsection{MRI contrast agent concentration - image signal relation}
\label{sec::singal-contrast}
Below, we briefly describe the relationship between the imaging signal 
seen in Figure~\ref{fig2} and the underlying MRI contrast agent
concentration. We remark that we use notations that are common in medical literature, which includes the use of two letter symbols. Hence, we will also use two letter symbols, such as $TE$ and $TR$, to keep the notation consistent with the presentation in~\cite{GOWLAND, MPRAGE}.   
The MRI contrast agent concentration $c$ causes the longitudinal(spin-lattice) relaxation time $T_{1}$ to shorten with the following relation
\begin{equation}
\frac{1}{T_{1}^{c}} = \frac{1}{T_{1}^{0}} + r_{1}c .
\label{EQ::contrast}
\end{equation}
The superscripts indicate relaxation time with MRI contrast agent $T_{1}^{c}$ and without MRI contrast agent $T_{1}^{0}$, while $r_1$ is the relaxivity constant for the MRI contrast agent in a medium. 
The MRI observations were obtained using a MRI sequence known as Magnetization Prepared Rapid Acquisition Gradient Echo (MPRAGE) with an inversion prepared magnetization. The relation between the signal and the relaxation time is non-linear, and is expressed with the following equations. The MRI signal value $S$ for this sequence can be expressed as
\begin{equation}
S = M_{n} \sin \theta e^{ - TE/T_2^{*} },
\label{EQ::SI_T2}
\end{equation}
with $TE$ and $\theta$ respectively denoting the echo time and the flip angle, and $M_{n}$ is the magnetization for the n-echo that we described below. 
Also $T_2^{*}$ is the transverse magnetization caused by a combination of spin-spin relaxation and magnetic field inhomogeneity, defined as 
\begin{equation}
\frac{1}{T_2^{*}} = \frac{1}{T_2} + \gamma \Delta B_{in} .
\end{equation}
Here $T_2$ is the transverse (spin-spin) relaxation time, $\gamma$ is the gyromagnetic ratio and $\Delta B_{in}$ is the magnetic field inhomogeneity across a voxel. The expression~\eqref{EQ::SI_T2} can be simplified by neglecting the exponential term, since $TE <<T_2^{*}$ is a general trait for this MRI sequence. Thus, \eqref{EQ::SI_T2} becomes 
\begin{equation}
S = M_{n} \sin \theta.
\label{EQ::SI}
\end{equation}
Magnetization for the n-echo $M_n$ is defined as~\cite{GOWLAND}:  
\begin{equation}
M_{n} = M_{0}  \left[ (1-\beta)\frac{(1-(\alpha \beta)^{n-1} }{1-\alpha\beta} + (\alpha \beta)^{n-1}(1-\gamma) + \gamma ( \alpha \beta)^{n-1} \frac{M_{e}}{M_{0}}  \right]   
\end{equation}
with 
\begin{equation}
\frac{M_{e}}{M_{0}} = - \left[ \frac{ 1 -\delta + \alpha \delta (1-\beta ) \frac{1-\alpha\beta^{m}}{1-\alpha \beta} + \alpha\delta(\alpha\beta)^{m-1} - \alpha^{m}\rho}{1 +\rho \alpha^{m} } \right].
\end{equation}
Using the following definitions
\begin{equation}
\begin{aligned}
\alpha &= \cos ( \theta ) \\
\beta  &= e^{- ^{T_b}/_{T_1^{c}} } \\
\delta &= e^{- ^{T_a}/_{T_1^{c}} } \\
\gamma &= e^{- ^{T_w}/_{T_1^{c}} } \\
\rho   &= e^{- ^{TR}/_{T_1^{c}}}  \\
T_w    &= TR - T_a -T_b(m-1)       .\\
\end{aligned}
\end{equation}
Here $T_b$ is known as the echo spacing time, $T_a$ is the inversion time, $T_w$ the time delay, $TR$ as the repetition time, $m$ is the number of echoes and $M_0$ is a calibration constant for the magnetization. The center echo denoted as $n=m/2$ will be the signal that we will consider when estimating concentration of the MRI contrast agent. Given~\eqref{EQ::SI}, the relative signal increase can be written as 
\begin{equation}
\frac{S^{c}}{S^{0}} = \frac{ M_{n}^{c} \sin (\theta)}{ M_{n}^{0} \sin (\theta) }.
\end{equation}
We define that  
\begin{equation}
f(T_1) = M_{n}/M_{0}. 
\label{scaledmagnetization}
\end{equation}
Figure~\ref{figuredti} shows $f(T_1)$ in CSF, grey and white matter. 
This gives the following relation 
\begin{equation}
\frac{f(T_{1}^{c} ) }{f(T_{1}^{0})}  = \frac{S^{c}}{S^{0}} 
\end{equation}
The signal difference between observation times were adjusted in~\cite{ringstad2018brain}. Thus we can express the change in $T_1$ due to MRI contrast agent as 
\begin{equation}
f ( T_{1}^{c} ) = \frac{S^{c}}{S^{0}} f(T_{1}^{0}) 
\end{equation}
and then estimate the concentration using~\eqref{EQ::contrast}. The $T_{1}^{0}$ values were obtained by T1 mapping of the brain using a MRI sequence known as MOLLI5(3)3~\cite{taylor2016t1}. This takes into account patient specific characteristic, such as tissue damage. Tissue damage can be observed in the MRI due to a lower signal in the white matter compared to healthy white matter tissue, thus damaged tissue have different $T_1$ relaxation time. 
The MRI contrast agent concentration was estimated in a preprocessing step, using the parameters obtained from the T1-map, MPRAGE MRI protocol~\cite{ringstad2018brain} and the value for $r_1$ found in~\cite{rohrer2005comparison}. The values of the function~\eqref{scaledmagnetization} was computed for $ T_1\in( 200, 4000)$ creating a lookup table. The lookup table was utilized with the baseline signal increase to estimate $T_1^{c}$, and then the concentration was computed using~\eqref{EQ::contrast}.  

\begin{figure}[h]
\centering
\includegraphics[width=0.70\textwidth]{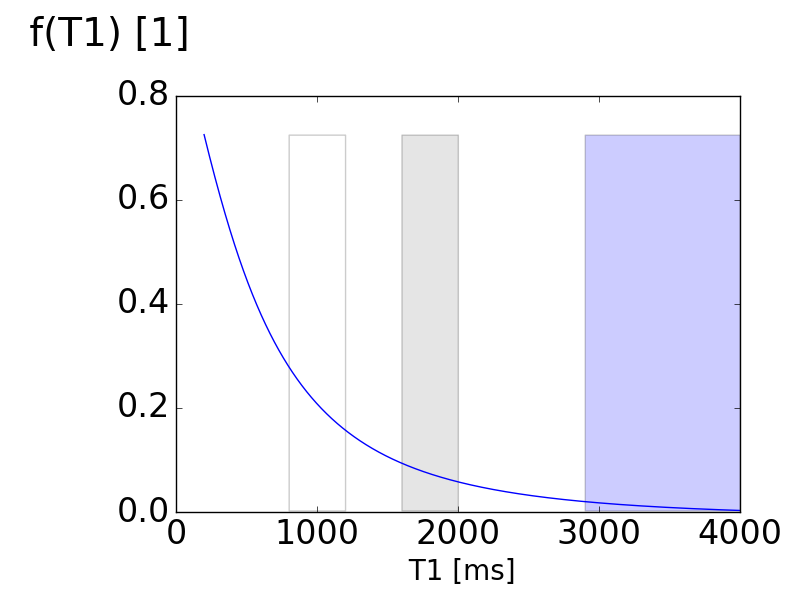} 
\caption{The image shows the function defined in~\eqref{scaledmagnetization} where the white region indicates $T_1$ values for white matter, the grey region indicates $T_1$ values for grey matter, the blue region indicates $T_1$ values for CSF.  }
\label{figureF} 
\end{figure}

\subsection{Diffusion tensor imaging}
\label{sec:supp:dti}
\begin{figure}
\centering
\includegraphics[width=0.75\textwidth]{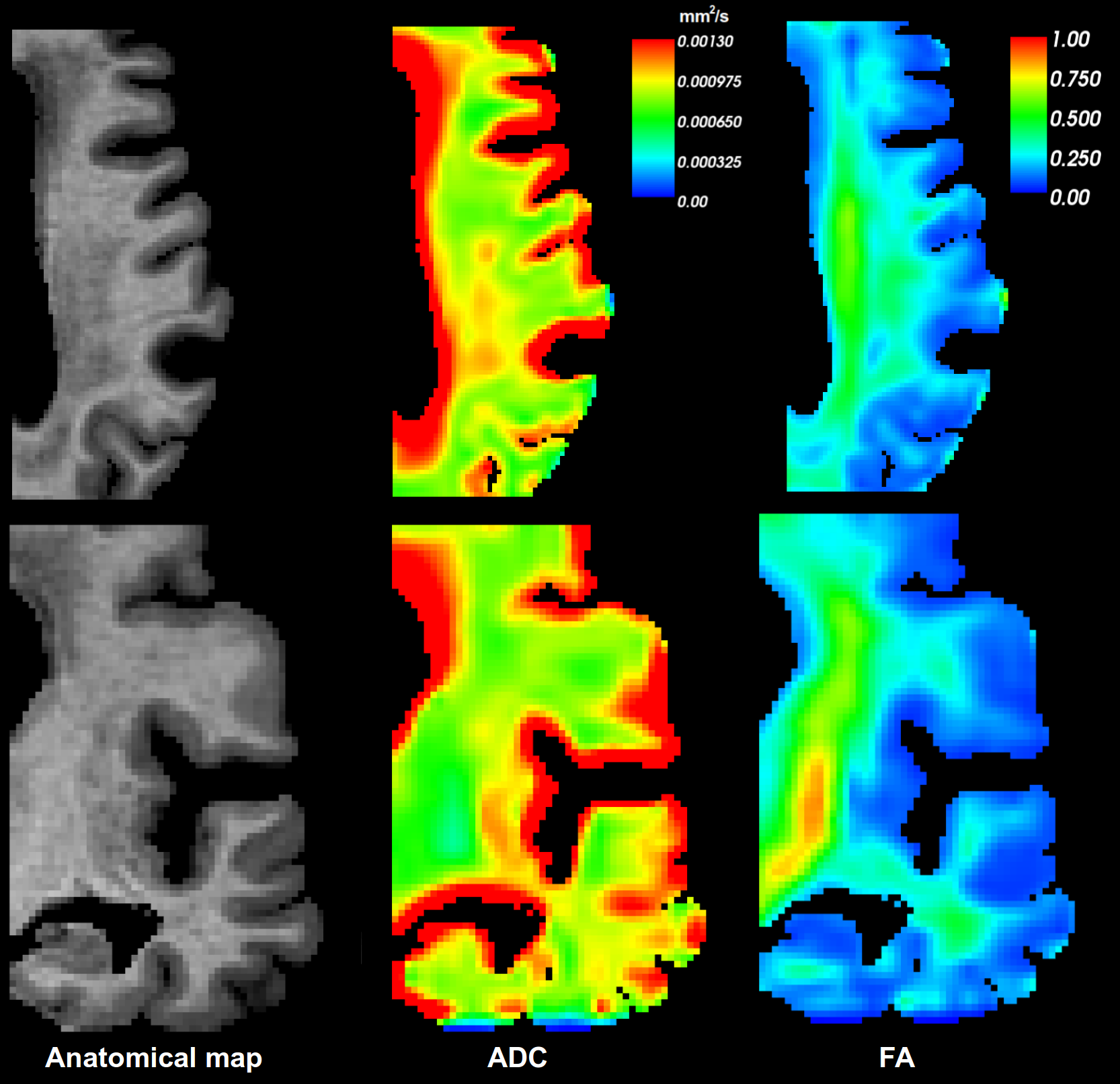} 
\caption{The left panel shows the anatomical map. The middle panel shows the apparent diffusion coefficients (ADC) obtained from DTI. The right panel shows the computed fractional anisotropy (FA) from the DTI at a color coded scale.}
\label{figuredti} 
\end{figure}

\begin{figure}
\centering
\includegraphics[scale=0.21]{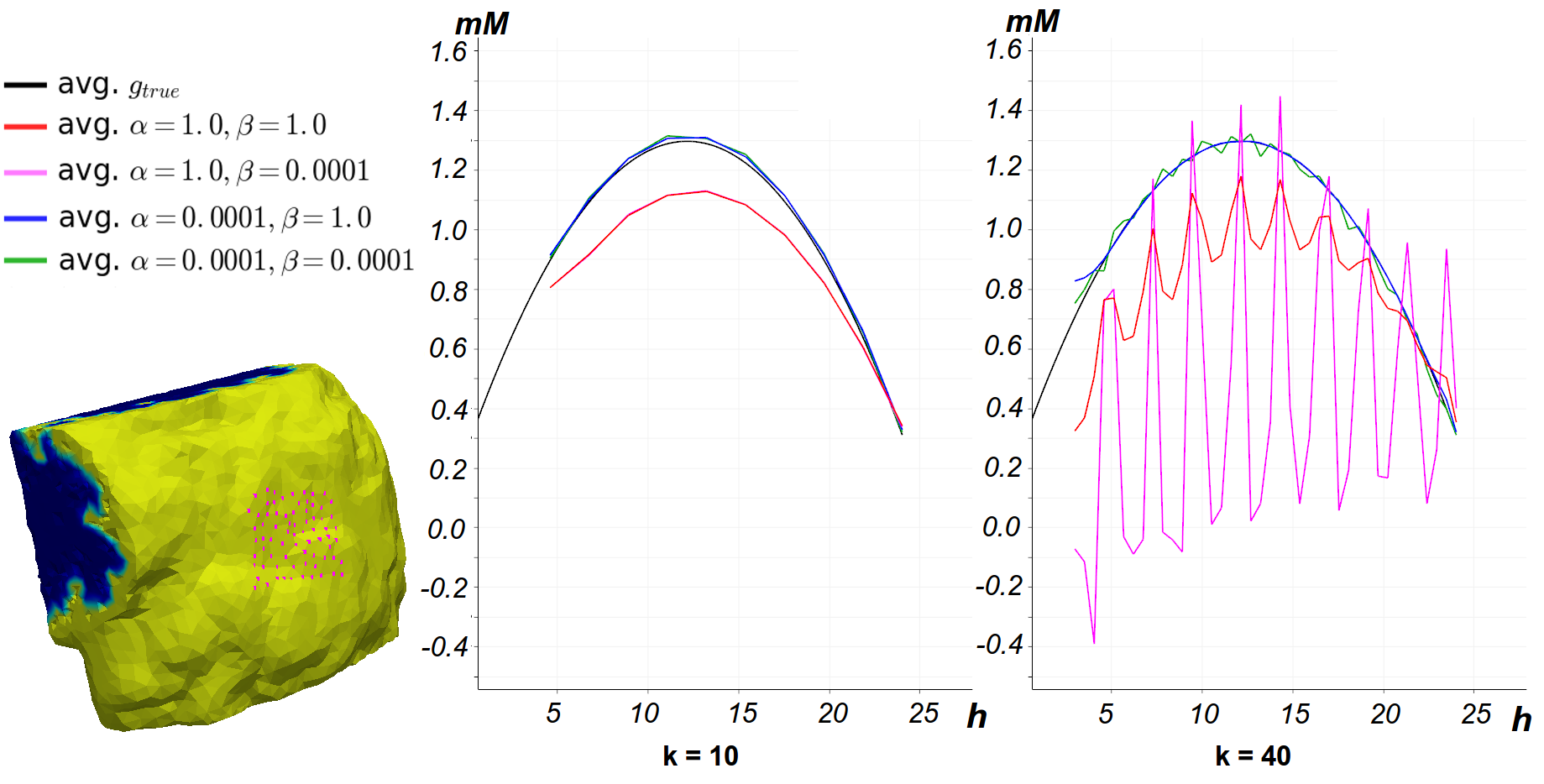}  
\caption{The image displays plots over time for a selection of points at the boundary of $\Omega_1$ with different regularization parameters and number of time steps $k$. The left panel shows the legend for the plot over time, together with the selection of points. The middle panel shows the average boundary value $g$ for different regularization parameter with $k=10$. The right panel shows the average boundary value $g$ for different regularization parameter with $k=40$.  }
\label{boundarycontrol}
\end{figure}

\begin{figure}
\centering
\includegraphics[scale=0.18]{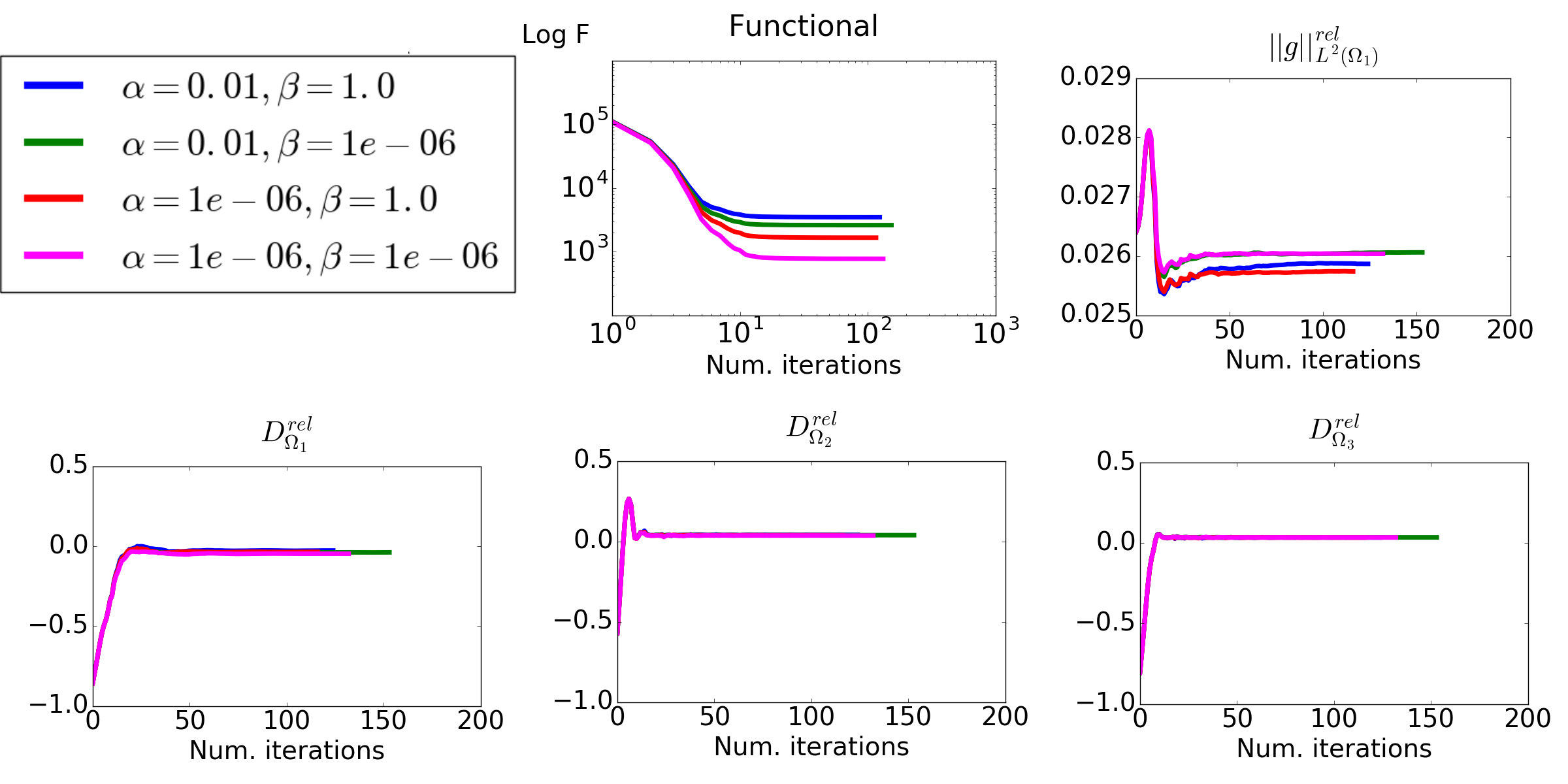}   
\caption{Convergence plots of the diffusion coefficients, boundary conditions and functional~\eqref{EQ::objf} with respect to different $\alpha$ and $\beta$ values. } 
\label{convergence}
\end{figure}

DTI images provide apparent diffusion coefficients (ADC) for water molecules (18Da) on short time-scales. A processed DTI image is shown in Figure~\ref{figuredti} with the largest ADC (shown in red in the middle figure) to be around $1.3\mathrm{e}{-3}  \mathrm{mm^2/s}$.  

In order to compare the computed diffusion coefficient, we need to estimate the ADC for Gadovist (604 Da)~\cite{MGadobutrol} by utilizing~\eqref{tortuosity} together with the tortuosity and the free diffusion coefficient for Gadovist. The tortuosity in the white and grey matter was estimated using~\eqref{tortuosity} with the ADC obtained from DTI and the self-diffusion of water which is $3.0\times 10^{-3}\mathrm{mm^{2}/s}$ at $37^{o}C$\cite{le2012diffusion}. We used the free diffusion coefficient for Gd-DPTA (550 Da)~\cite{MGgDPTA}, which was reported in~\cite{GdDPTA-DIFFUSION} to be $3.8\times 10^{-4} \mathrm{mm^{2}/s}$, as a surrogate for Gadovist. This is based on the Stokes-Einstein equation, which states that molecules with similar weight will have similar diffusion coefficients~\cite{valencia2011understanding}. Furthermore, we need to assume that the tortuosity is independent for molecules with mass lower than 1kDa, i.e. not large proteins.
The fractional anisotropy is defined as 
\begin{equation}
FA^{2} =  \frac{3}{2} \frac{ (\lambda_1 - MD )^{2} +(\lambda_2 - MD )^{2} +(\lambda_3 - MD )^{2}}{\lambda^{2}_1 + \lambda^{2}_2  +\lambda^{2}_3 },
\end{equation}
with the mean diffusivity $MD$, closely related to ADC, defined as 
\begin{equation}
MD = \frac{\lambda_1 +\lambda_2 +\lambda_3 }{3}.
\end{equation}
In these equations $\lambda_i$ denotes the eigenvalues of the diffusion tensor.
$FA$ is a common measure of anisotropy and shown in the right-most image in Figure~\ref{figuredti} for our patient.

\end{document}